\documentstyle [12pt,leqno]{amsart}

\def\O{\mbox{${\cal O}$}}
\def\M{\mbox{${\cal M}$}}
\def\B{\mbox{${\cal B}$}}
\def\T{\mbox{${\rm Def}(X)$}}
\def\P{\mbox{${\Bbb P}$}}
\def\l{\lambda}
\def\a{\alpha}
\def\w{\omega}

\title{Rational Curves on Calabi-Yau Threefolds}
\author{Sheldon Katz}
\address{Department of Mathematics, Oklahoma State University, Stillwater, OK
74078}
\thanks{On leave 1991--92 at  Department of Mathematics,
Duke University, Durham, NC 27706}

\begin{document}

\maketitle

In the conformal field theory arising from the compactification of
strings on a Calabi-Yau threefold $X$, there naturally arise
fields correspond to harmonic forms of types (2,1) and (1,1) on $X$
\cite{CHSW}.
The uncorrected Yukawa couplings in $H^{2,1}$ and $H^{1,1}$ are cubic forms
that can be constructed by techniques of algebraic geometry --- \cite{Dave}
contains a nice survey of this in a general context in a language written for
mathematicians.  The cubic form on the space $H^{p,1}$ of harmonic $(p,1)$
forms
is given by the intersection product $(\w_i,\w_j,\w_k)\mapsto\int_X\w_i\wedge
\w_j\wedge\w_k$ for $p=1$, while for $p=2$ there is a natural formulation in
terms of infinitesimal variation of Hodge structure \cite{Dave,IVHS}.  The
Yukawa couplings on $H^{2,1}$ are exact, while those on $H^{1,1}$ receive
instanton corrections.  In this context, the {\em instantons} are non-constant
holomorphic maps $f:{\Bbb C}{\Bbb P}^1\to X$.  The image of such a map is a
{\em rational curve} on $X$, which may or may not be smooth.  If the rational
curve $C$ does not move inside $X$, then the contribution of the instantons
which
are generically 1-1 (i.e. birational) maps with image $C$ can be written down
explicitly --- this contributions only depends on the topological type of $C$,
or more or less equivalently, the integrals $\int_Cf^*J$ for (1files,
a,1) forms $J$
on
$X$.

If the conformal field theory could also be expressed in terms of a ``mirror
manifold'' $X'$, then the uncorrected Yukawa couplings on $H^{2,1}(X')$
would be the same as the corrected Yukawa couplings on $H^{1,1}(X)$.  So if
identifications could be made properly, the infinitesimal variation of Hodge
structure on $X'$ would give information on the rational curves on $X$.  In a
spectacular paper \cite{Ccalc}, Candelas et.\ al.\ do this when $X$ is a
quintic threefold.  Calculating the Yukawa coupling on $H^{2,1}(X')$ as the
complex structure of $X'$ varies gives the Yukawa coupling on $H^{1,1}(X)$
as $J$ varies.  The coefficients of the resulting Fourier series are then
directly related to the instanton corrections.

To the mathematician, there are some unanswered questions in deducing the
number
of rational curves of degree $d$ on $X$ from this \cite{Dave}.  I merely cite
one problem here, while I will state a little differently.  It is not yet known
how to carry out the calculation of the instanton
correction associated to a continuous family of rational curves on $X$.
If there were continuous families of rational curves on a general quintic
threefold $X$,
this would complicate
the instanton corrections to the Yukawa coupling on the space of $(1,1)$ forms.
The question of whether or not such families exist has not yet been resolved.
But even worse, this complication is present in any case --- any degree $m$
mapping of
${\Bbb C}{\Bbb P}^1$ to itself may be composed with any birational mapping of
${\Bbb C}{\Bbb P}^1$ to $X$ to give a new map from ${\Bbb C}{\Bbb P}^1$ to $X$;
and this
family has
has moduli than mere reparametrizations of the sphere.  In \cite{Ccalc} this
was
explicitly recognized; the
assertion was made that such a family counts $1/{m^3}$ times.  This has since
been verified by Aspinwall and Morrison \cite{TFT}.

My
motivation in writing this note is to give a general feel for the
mathematical meaning of ``the number of rational curves on a Calabi-Yau
threefold'', in particular, how to ``count'' a family of rational curves.
I expect that these notions will directly correspond to the not as yet worked
out procedure for calculating instanton corrections associated to general
families of rational curves.  More satisfactory mathematical formulations are
the subject of work still in progress.

I want to take this opportunity to thank David Morrison for his suggestions, to
the Mathematics Department at Duke University for its
hospitality while this manuscript
was being prepared, and to the organizers of the Mirror Symmetry Workshop for
providing a wonderful opportunity for mathematicians and physicists to learn
from each other.

\section{Rational Curves, Normal Bundles, Deformations}

Consider a Calabi-Yau threefold $X$ containing a smooth rational curve
\linebreak
$C\cong{\Bbb C}{\Bbb P}^1$.  The normal bundle $N_{C/X}$ of $C$ in $X$ is
defined by the
exact sequence

\begin{equation}
\label{nbseq}
0\to T_C\to T_X|_C\to N_{C/X}\to 0.
\end{equation}

$N=N_{C/X}$ is a rank 2 vector bundle on $C$, so $N=\O(a)\oplus\O(b)$ for some
integers $a,b$.  Now $c_1(T_C)=2$, and $c_1(T_X)=0$ by the Calabi-Yau
condition.
So the exact sequence (\ref{nbseq}) yields $c_1(N)=-2$, or $a+b=-2$.  One
``expects'' $a=b=-1$ in the general case.  This is because there is a moduli
space of deformations of the vector bundle $\O(a)\oplus\O(b)$, and the general
point of this moduli space is $\O(-1)\oplus\O(-1)$, no matter what $a$ and $b$
are, as long as $a+b=-2$.

Let \M\ be the moduli space of rational curves in $X$.  The tangent space to
\M\ at
$C$ is given by $H^0(N)$ \cite[\S 12]{KS}.  In other words, \M\ may be locally
defined by finitely many equations in $\dim H^0(N)$ variables.

\bigskip\noindent
{\bf Definition.}
$C$ is {\em infinitesimally rigid} if $H^0(N)=0$.

\bigskip
Infinitesimal rigidity means that $C$ does not deform inside $X$, not even to
first order.
Note that $H^0(N)=0$ if and only if $a=b=-1$.  Thus

\begin{itemize}
\item $C$ is infinitesimally rigid if and only if $a=b=-1$.
\item $C$ deforms, at least infinitesimally, if and only if $(a,b)\neq
(-1,-1)$.
\end{itemize}

\M\ can split up into countably many irreducible components.  For instance,
curves
with distinct homology classes in $X$ will lie in different components of \M.
However, there will be at most finitely many components of \M\ corresponding to
rational curves in a fixed homology class.

There certainly exist Calabi-Yau threefolds $X$ containing
positive dimensional families of rational
curves.  For instance, the Fermat quintic threefold
$x_0^5+\ldots x_4^5=0$ contains the family of lines given parametrically in the
homogeneous coordinates $(u,v)$ of ${\Bbb P}^1$ by
$(u,-u,av,bv,cv)$, where $(a,b,c)$ are the parameters of the plane curve
$a^5+b^5+c^5=0$.  However, suppose that all rational curves on $X$ have normal
bundle $\O(-1)\oplus\O(-1)$.  Then since \M\ consists entirely of
discrete points, the remarks above show that there would be only finitely many
rational
curves in $X$ in a fixed homology class.  Since a quintic threefold has
$H_2(X,{\Bbb Z})\cong{\Bbb Z}$, the degree of a curve is essentially the same
as its
homology class.  This discussion leads to Clemens' conjecture \cite{Ccon}:

\bigskip\noindent
{\bf Conjecture (Clemens).}
{\em A general quintic threefold contains only finitely many rational curves of
degree
$d$, for any $d\in {\Bbb Z}$. These curves are all infinitesimally rigid.}

\bigskip
Clemens' original constant count went as follows \cite{Chom}:
a rational curve of degree $d$ in
$\P^4$ is given parametrically by 5 forms  $\a_0(u,v),\ldots,\a_4(u,v)$, each
homogeneous of
degree $d$ in the homogeneous coordinates $(u:v)$ of $\P^1$.  These $\a_i$
depend on $5(d+1)$
parameters.  On the
other hand, a quintic equation $F(x_0,\ldots,x_4)=0$ imposes the condition
$F(\a_0(u,v),\ldots,\a_4(u,v))\equiv 0$ for the parametric curve to be
contained in this
quintic threefold.  This is a polynomial equation of degree $5d$ in $u$ and
$v$.
Since a general degree $5d$ polynomial $\sum a_{i}u^iv^{5d-i}$ has $5d+1$
coefficients, setting these equal to zero results in $5d+1$ equations among the
$5d+5$ parameters of the $\a_i$.  If $F$ is {\em general}, it seems plausible
that
these equations should impose independent conditions, so that the solutions
should
depend on $5d+5-(5d+1)=4$ parameters.  However, any curve has a 4-parameter
family of reparametrizations $(u,v)\mapsto (au+bv,cu+dv)$, so there are
actually
a zero dimensional, or finite number, of curves of the general $F=0$.

The conjecture is known to be true for $d\le 7$ \cite{Kfin}.  For any $d$, it
can even be proven that there exists an infinitesimally rigid curve of degree
$d$ on a general $X$.  Similar conjectures
can be stated for other Calabi-Yau threefolds.

There are many kinds of non-rational
curves which appear to occur in finite number on a general quintic
threefold.  For instance,
elliptic cubic curves are all planar.  The plane $P$ that one spans meets the
quintic in a quintic curve containing the cubic curve.  The other component
must
be a conic curve.  This sets up a 1-1 correspondence between elliptic cubics
and
conics on any quintic threefold.  Hence the number of elliptic cubics on a
general
quintic must be the same as the number of conics, 609250 \cite{Kfin}.
Finiteness of
elliptic quartic curves has been proven by Vainsencher \cite{V}; the actual
number has not yet been computed.

On the other hand, there are infinitely many plane quartics on {\em any\/}
quintic
threefold: take any line in the quintic, and each of the infinitely many planes
containing the line must meet the quintic in the original line union a quartic.

If a curve has $N\cong\O\oplus\O(-2)$, then $C\subset X$ deforms to first
order.  In fact, \nolinebreak since $H^0(N)$ is one-dimensional, there is a
family of curves on $X$ parametrized by a single variable $t$,
subject to the constraint $t^2=0$.  In other words, start with a rational curve
given parametrically by forms $\a_0,\ldots ,\a_4$, homogeneous of degree $d$
in $u$ and $v$.  Take a perturbation
$\a_i(u,v;t)=\a_i(u,v)+t\a_i'(u,v)$, still homogeneous in $(u,v)$.
Form the equation
$F(\a_0,\ldots a_4)=0$ and formally set $t^2=0$; the resulting equation has a
5 dimensional space of solutions for the $\a_i'$, which translates into a
unique
solution up to multiples and reparametrizations of $\P^1$.  The curve $C$
but may
or may not deform to second order.  $C$ deforms to $n^{\mbox{th}}$ order for
all $n$ if and only if $C$ moves in a 1-parameter family.  A pretty description
of the general situation is given in \cite{R}.

If $C$ deforms to
$n^{\mbox{th}}$ order, but not to $(n+1)^{\mbox{th}}$ order, then  one sees
that
while $C$ is an isolated point in the moduli space of curves on $X$, it
more naturally is viewed as the solution to the equation $t^{n+1}=0$ in one
variable $t$.  So $C$
should be viewed as a rational curve on $X$ with multiplicity $n+1$.

If a curve has $N\cong\O(1)\oplus\O(-3)$, $C$ has a 2 parameter space of
infinitesimal deformations, and the structure of \M\ at $C$ is correspondingly
more
complicated.  An example is given in the next section.  The general situation
has
not yet been worked out.

\section{Counting Rational Curves}

In this section, a general procedure for calculating the number of smooth
rational curves of a
given type is described. Alternatively, a canonical definition of
this number can be given using the Hilbert scheme (this is what was done
by Ellingsrud and Str{\o}mme in their work on twisted cubics \cite{EScalc});
however, it is usually quite difficult to implement a calculation along
these lines.

Embed the Calabi-Yau threefold $X$ in a larger compact space
${\Bbb P}$ (which may be thought of as a projective space, a weighted
projective
space, or a product of such spaces).  $\M_\l$ will denote the moduli space
parametrizing smooth rational curves in ${\Bbb P}$ of a given topological type
or degree $\lambda$.
\T\ denotes the irreducible component of $X$ in the moduli space of Calabi-Yau
manifolds
in ${\Bbb P}$ (here the K\"ahler structure is ignored).  In other words, \T\
parametrizes the deformations of $X$ in ${\Bbb P}$.

\begin{enumerate}

\item Find a compact moduli space $\bar {\M_\l}$ containing $\M_\l$ as a
dense open subset,
such that the points of $\bar{\M_\l}-\M_\l$ correspond to degenerate
curves of type $\l$. $\bar {\M_\l}$ parametrizes degenerate deformations of
the smooth curve (not the mapping from ${\Bbb C}{\Bbb P}^1$ to $X$).
It is better for $\bar{\M_\l}$ to be smooth.

\pagebreak
\item Find a rank $r=\dim(\M_\l)$ vector bundle $\B$ on
$\bar{\M_\l}$ such that
  \begin{enumerate}
  \item To each $X'\in\T$ there is a section $s_{X'}$ of \B\ which vanishes at
  	$C\in\bar{\M_\l}$ if and only if $C\subset X'$.
  \item There exists an $X'\in\T$ such that $s_{X'}(C)=0$ if and only if
  	$C\in \M_\l$ and $C\subset X'$.
 \item $C$ is an isolated zero of $s_{X'}$ with $\mbox{mult}_C(s_{X'})=1$ if
	and only if
  	$N_{C/X}\cong\O(-1)\oplus\O(-1)$.
  \end{enumerate}
\end{enumerate}

\bigskip\noindent
{\bf Working Definition.} The {\em number of smooth rational curves $n_\l$ of
type $\l$}
is given by the
$r^{\rm th}$ Chern class $c_r(\B)$.

\bigskip

Why is this a reasonable definition?
Suppose that an $X'\subset\T$ can be found with the properties required above,
with the additional property that there are only finitely many curves of type
$\l$ on $X'$,
and that they are all infinitesimally rigid.  Then it can be checked that the
number of
(possibly degenerate) curves
of type $\l$ on $X'$ is independent of the choice of $X'$ satisfying the above
properties, and is also equal to $c_r(\B)$.  This last follows since
$c_{\mbox{rank}(E)}(E)$ always gives the homology class of the zero locus $Z$
of any
section of any bundle $E$ on any variety $Y$, whenever
$\dim(Z)=\dim(Y)-\mbox{rank}
(E)$.  In our case, $0=\dim(\{{\rm lines}\})=\dim(\M_\l)-\mbox{rank}(\B)$.  In
essentially all known cases, the number of curves has been worked out by the
method of this working definition.  Examples are given below.

There is a potential problem with this working definition.  For
families of Calabi-Yau threefolds such that no threefold in the family contains
finitely many curves of given type, it may be that the ``definition'' depends
on the choice
of compactification and/or vector bundle, i.e. this is not well-defined.  My
reason for almost calling this method a definition is that it {\em does\/}
give a finite number corresponding to an infinite family of curves, which
{\em is\/}
well-defined in the case that the Calabi-Yau threefold in question belongs to
a family containing some other Calabi-Yau threefold with only finitely many
rational curves of the type under consideration.

In the case where the general $X$ contains irreducible singular curves which
are
the images of maps from $\P^1$, a separate
but similar procedure must be implemented to calculate these, since they give
rise to instanton corrections as well.  For example,
there is a 6-parameter family of two-planes in $\P^4$.  A two plane $P$ meets a
quintic threefold $X$ in a plane quintic curve.  For general $P$ and $X$,
this curve is a smooth genus 6 curve.  But if the curve acquired 6 nodes, the
curve would be rational.  This being 6 conditions on a 6 parameter family, one
expects that a general $X$ would contain finitely many 6-nodal rational plane
quintic curves, and it can be verified that this is indeed the case.
The problem is
that there is no way to deform a smooth rational curve to such a singular
curve --- the dimension of $H^1(\O_C)$ is zero for a smooth rational curve $C$,
but positive for such singular curves, and this dimension is a deformation
invariant \cite[Theorem 9.9]{H}.  For example, if one tried to deform a smooth
twisted cubic curve to a
singular cubic plane curve by projecting onto a plane, there would result an
``embedded point'' at the singularity, creating a sort of discontinuity in the
deformation process \cite[Ex.\ 9.8.4]{H}.

\bigskip\noindent
{\em Examples:}
\begin{enumerate}
\item Let $X\subset\P^4$ be a quintic threefold.  Take $\l=1$, so that we are
counting lines.  Here $\M_1=G(1,4)$ is the Grassmannian of lines in $\P^4$ and
is
already compact, so take $\bar{\M_1}=\M_1$.  Let $\B=\mbox{Sym}^5(U^*)$, where
$U$, the universal bundle, is the rank 2 bundle on $\M_1$ whose fiber over a
line
$L$ is the 2-dimensional subspace $V\subset{\Bbb C}^5$ yielding $L\subset\P^4$
after
projectivization.  Note that $\mbox{rank}(\B)=\dim(\M_1)=6$.
$\T\subset\P(H^0(\O_{\P^4}(5)))$ is the subset of smooth
quintics.  A quintic $X$ induces a section $s_X$ of \B, since an equation for
$X$ is a
quintic form on ${\Bbb C}^5$, hence induces
a quintic form on $V$ for $V\subset{\Bbb C}^5$ corresponding to $L$.  Clearly
$s_X(L)=0$ if and only if $L\subset X$.  The above conditions are easily seen
to
hold.   $c_6(\B)=2875$ is the number of lines on $X$.  This
calculation is essentially the same as that done for cubics in
\cite[Thm.\ 1.3]{AK},
where dual notation is used, so that the $U^*$ used here becomes the universal
quotient bundle $Q$ in \cite{AK}.  The number $2875$ agrees with the result of
Candelas et.\ al.\ \nopagebreak\cite{Ccalc}.

\bigskip
\item Continuing with the quintic, take $\l=2$.  Any conic $C$ necessarily
spans a
unique 2-plane containing $C$.  Let $G=G(2,4)$ be the Grassmannian of 2-planes
in
$\P^4$, and let $U$ be the rank 3 universal bundle on $G$.  Put
$\bar{\M_2}=\P(\mbox{Sym}^2(U^*))$ be the projective bundle over $G$ whose
fiber
over a plane $P$ is the projective space of conics in $P$.  Clearly
$\M_2\subset\bar{\M_2}$ (but they are not equal; $\bar{\M_2}$ also contains the
union of any two lines or a double line in any plane).  Let
$\B=\mbox{Sym}^5(U^*)/(\mbox{Sym}^3(U^*)\otimes\O_{\P}(-1))$ be the bundle on
$\bar{\M_2}$ of quintic
forms on the 3 dimensional vector space $V\subset{\Bbb C}^5$, modulo those
which
factor as any cubic times the given conic. ($\O_{\P}(-1)$ is the line bundle
whose
fiber over a conic is the one dimensional vector space of equations for the
conic within its supporting
plane.  The quotient is relative to the natural embedding
$\mbox{Sym}^3(U^*)\otimes\O_{\P}(-1)\to\mbox{Sym}^5(U^*)$ induced by
multiplication.) $\mbox{rank}(\B)=\dim(\bar{\M_2})=11$.
$c_{11}(\B)=609250$.  See \cite{Kfin} for more details, or \cite{CMW} for an
analogous computation in the case of quartics.  The number $609250$ agrees
with
the result of
Candelas et.\ al.

\bigskip
\item Again consider the quintic, this time with $l=3$.  In \cite{EScalc},
Ellingsrud and Str{\o}mme take $\bar{\M_3}$ to be the closure of the locus of
smooth twisted cubics in $\P^4$ inside
the Hilbert scheme.  This space has dimension 16.  $\B$ is essentially the rank
16 bundle of degree 15 forms on $\P^1$ induced from quintic polynomials in
$P^4$
by the the degree 3 parametrization of the cubic (it must be shown that this
makes sense for
degenerate twisted cubics as well). The equation of a general quintic gives a
section of $\B$, vanishing precisely on the set of cubics contained in $X$.
$c_{16}(\B)=317206375$, again agreeing with
the result of
Candelas et.\ al.

\end{enumerate}

\bigskip

The key to the first two calculations are the Schubert calculus for calculating
in Grassmannians \cite[Ch.\ 1.5]{GH} and standard formulas for projective
bundles \cite[Appendix A.3]{H}.  The third calculation is more intricate.

Regarding complete intersection Calabi-Yau manifolds, similar examples are
found in
\cite{Klin} for lines  and \cite{SvS} for conics.

Note that the number of curves $c_r(\B)$ in no way depends on the choice
of $X$, even if $X$ is a degenerate Calabi-Yau threefold, or contains
infinitely many rational curves of type $\l$.  It turns out that a natural
meaning can be assigned to this number.  In fact, the moduli space of
curves of type $\l$ on $X$ splits up into ``distinguished varieties''
$Z_i$
\cite{F}, and a number, the {\em equivalence} of $Z_i$, can be assigned
to each distinguished variety (the number is 1 for each
$\O(-1)\oplus
\O(-1)$ curve).   This number is precisely equal to the number of curves on $X$
in $Z_i$
which arise as limits of curves on $X'$ as $X'$ approaches $X$ in a 1-parameter
family \cite[Ch.\ 11]{F}.

\bigskip\noindent
{\em Examples:}
\begin{enumerate}
\item If a quintic threefold is a union of a hyperplane and a quartic, then the
quintic contains infinitely many lines and conics.  However, in the case of
lines,
given a general 1-parameter family of quintics approaching this reducible
quintic,
$1275$ lines approach the hyperplane, and 1600 lines approach the quartic
\cite{Kdeg}.
So the infinite set of lines in the hyperplane ``count'' as $1275$, while
those in
the quartic count as $1600$.  For counting conics, the $609250$ conics
distribute
themselves as $187250$ corresponding to the component of conics in the
hyperplane,
$258200$ corresponding to the component of conics in the quartic, and $163200$
corresponding to the component of conics which degenerate into a line in the
hyperplane union an intersecting line in the quartic \cite{Kit}.  Note that
this
reducible conic lies in $\bar{\M_2}-\M_2$, and illustrates why $\M_2$ itself is
insufficient for calculating numbers when there are infinitely many curves.
Most of these numbers have been calculated recently by Xian Wu \cite{W} using
a different method.

\item If a quintic threefold is a union of a quadric and a cubic, the lines
on the quadric count as $1300$, and the lines on the cubic count as $1575$
\cite{Kdeg,W2}.  The conics on the quadric count as $215,950$, while the conics
on
the cubics count as $243900$ \cite{W2}. Presumably this implies that the conics
which
degenerate into a line in the quadric union an intersecting line in the cubic
count as $609250-(215950+243900)=149400$, but this has not been checked
directly yet.

\item There are infinitely many lines on the Fermat quintic threefold
$x_0^5+\ldots+x_4^5=0$.  These divide up into 50 cones, a typical one being
the family of lines given parametrically in the
homogeneous coordinates $(u,v)$ of ${\Bbb P}^1$ by
$(u,-u,av,bv,cv)$, where $(a,b,c)$ satisfy $a^5+b^5+c^5=0$. Each of these count
as 20.  There are also 375 special lines, a typical one being given by the
equations
$x_0+x_1=x_2+x_3=x_4=0$ (these lines were also noticed in \cite{DSWW}).  These
lines $L$ count with multiplicity 5.  Note that $50\cdot 20+375\cdot 5=
2875$ \cite{AKlin}.  This example illustrates the potential complexity in
calculating the distinguished varieties
$Z_i$ --- some components can be embedded inside others.  This may be
understood as well by looking at the moduli space of lines on the Fermat
quintic locally
at a line corresponding to a special line.
Since \linebreak $N_{L/X}\cong\O(1)\oplus\O(-3)$, $H^0(N_{L/X})$ has dimension
2, and the
moduli space of lines is a subset of a 2-dimensional space.
A calculation shows that it
is locally defined inside 2
dimensional $(x,y)$ space by the equations $x^2y^3=x^3y^2=0$.  The $x$ and $y$
axes
correspond to lines on each of the 2 cones, each occurring with multiplicity
2; there
would be just one equation $x^2y^2=0$ if the special line corresponding to
$(0,0)$ played
no role; since this is not the case, it can be expected to have its own
contribution; i.\hskip0pt e.\ each special line is a distinguished
variety.

\item Examples for lines in complete intersection Calabi-Yau threefolds were
worked
out in \cite{Klin}.

\end{enumerate}

\bigskip
These numbers can also be calculated by intersection
theory techniques.  Let  $s(Z_i,\bar{\M_\l})$ denote the {\em Segre class} of
$Z_i$
in $\bar{\M_\l}$.  If $Z_i$ is smooth, this is simply the formal inverse of the
total Chern class $1+c_1(N)+c_2(N)+\ldots$ of the normal bundle $N$ of $Z_i$
in $\bar{\M_\l}$.  Then if $Z_i$ is a connected component of the zero locus of
$s_X$, the equivalence of $Z_i$ is
the zero dimensional part of $c(\B)\cap s(Z_i,\bar{\M_\l})$
\cite[Prop.\ 9.1.1]{F}.

If $Z_i$ is an irreducible component which is not a connected component, then
this
formula is no longer applicable.  However, I have had recent success with a new
method that supplies ``correction terms'' to this formula.  The method is
currently ad hoc (the most relevant success I have had is in calculating the
``number'' of
lines on a cubic surface which is a union of three planes), but I expect that
a more
systematic procedure can be developed.

This of course is a reflection of the situation in calculating instanton
corrections to the Yukawa couplings.  If there is a continuous family of
instantons, then calculating the corrections will be more difficult.  If the
parameter space for instantons is smooth, this should make the
difficulties more managable.  If the
space is singular, the calculation is more difficult.  I expect that the
calculation would be even more difficult if instantons occur in at least 2
families that intersect.

Of course, in the calculation of the Yukawa couplings via path integrals, there
is no mention of vector bundles on the moduli space.  This indicates to
me that there should be a mathematical definition of the
equivalence of a distinguished variety that does not refer to an auxiliary
bundle.  Along these lines, one theorem will
be stated without proof.

Let $Z$ be a $k$-dimensional unobstructed family of rational curves on a
Calabi-Yau threefold
$X$.  There is the total space ${\cal Z}\subset Z\times X$ of the family, with
projection map $\pi:{\cal Z}\to Z$ such that $\pi^{-1}(z)$ is the curve in $X$
corresponding to $z$, for each $z\in Z$.  Let $N$ be the
normal bundle of ${\cal Z}$ in $Z\times X$.  Define the equivalence $e(Z)$ of
$Z$ to be the number $c_k(R^1\pi_*N)$.  For example, if $k=0$, then the curve
is infinitesimally rigid, and $e(Z)=1$.

\bigskip\noindent
{\bf Theorem.} {\em Let $Z$ be an unobstructed family of rational curves of
type $\l$ on a
Calabi-Yau threefold $X$.  Suppose that $X$ deforms to a Calabi-Yau threefold
containing only finitely many curves of type $\l$.  Then precisely $e(Z)$ of
these curves (including multiplicity) approach curves of $Z$ as the Calabi-Yau
deforms to $X$}.

\bigskip
Some of the examples given earlier in this section can be redone via this
theorem.  Also, as anticipated by \cite{Ccalc}, it can be calculated that a
factor of $1/{m^3}$ is introduced
by degree $m$ covers by a calculation similar in spirit to that found in
\cite{TFT}.

\end{document}